\input amstex
\documentstyle{amsppt}
\magnification=1200

\hoffset=-0.5pc
\vsize=57.2truepc
\hsize=38truepc
\nologo
\spaceskip=.5em plus.25em minus.20em
\define\cbergthe {1}
\define\bottone {2}
\define\botshust {3}
\define\curtione {4}
\define\dwykanon {5}
\define\eilmactw {6}
\define\eilmacth {7}
\define\gugenhtw {8}
\define\gugenmay {9}
 \define\huebkade{10}
 \define\lattice {11}
\define\kanone {12}
\define\maclafou {13}
\define\maclbotw {14}
\define\milgrone {15}
\define\milnoru {16}
\define\milnothr {17}
\define\mooreone {18}
\define\puppeone {19}
\define\quillone {20}
\define\gsegatwo {21}
\define\gsegathr {22}
\define\stashone {23}
\define\stasheig {24}
\define\steenone {25}
\define\weingone {26}
\define\wongone {27}
\topmatter
\title 
Comparison of the geometric bar and W-constructions\\
J. Pure and Applied Algebra 131 (1998) 109--123
\endtitle
\author Clemens Berger and Johannes Huebschmann{\dag}
\endauthor
\affil 
Universit\'e de Nice-Sophia Antipolis
\\
Laboratoire J. A. Dieudonn\'e, URA 168
\\
Parc Valrose
\\
06108 NICE Cedex 2
\\
France
\\
cberger\@math.unice.fr
\\
and
\\
Max Planck Institut f\"ur Mathematik
\\
Gottfried Claren Str. 26
\\
53225 BONN
\\
Germany
\\
huebschm\@mpim-bonn.mpg.de 
\endaffil
\thanks
{{\dag} The second named author carried out this work in the framework of the
VBAC research group of EUROPROJ.}
\endthanks
\abstract{For a simplicial group $K$, the
realization of the
$W$-construction
\linebreak
$WK \to \overline WK$ 
of $K$ is shown to be naturally homeomorphic
to the universal bundle 
$E|K| \to B|K|$ 
of its geometric realization $|K|$.
The argument involves certain recursive descriptions
of the $W$-construction and classifying bundle
and relies on the facts that the realization functor
carries an action of a simplicial group to a geometric
action of its realization and preserves
reduced cones and colimits.}
\endabstract
\address
{
\newline
Permanent address of the second named author:
\newline
Universit\'e de Lille 1, UFR de Math\'ematiques
\newline
UMR CNRS 8524, Laboratoire Paul Painlev\'e
\newline
Labex CEMPI (ANR-11-LABX-0007-01)
\newline
59655 Villeneuve d'Ascq Cedex/France
\newline
Johannes.Huebschmann\@math.univ-lille1.fr
}
\endaddress
\date{April 23, 1996} \enddate
\endkeywords

\endtopmatter
\document
\leftheadtext{Clemens Berger and Johannes Huebschmann}
\rightheadtext{Comparison of the geometric bar and W-constructions}
\beginsection Introduction

Let $K$ be a simplicial group;
its realization $|K|$
is a topological group
suitably interpreted
when $K$ is not countable.
The $W$-construction
$WK \to \overline WK$ yields a functorial universal
simplicial principal $K$-bundle,
and the classifying bundle construction
$E|K| \to B|K|$ 
of its geometric realization $|K|$
yields a
functorial universal
principal $|K|$-bundle.
The realization
of the $W$-construction
also yields a
universal
principal $|K|$-bundle
$|WK| \to |\overline WK|$,
by virtue of the general realization result in \cite\weingone.
In this note we identify the 
the classifying bundle
$E|K| \to B|K|$
with the realization of the $W$-construction.
A cryptic remark about the possible coincidence
of the two constructions may be found in the introduction
to Steenrod's paper \cite\steenone\ 
but to our knowledge this has never been made explicit
in the literature.
\smallskip
Spaces are assumed compactly generated,
and all constructions on spaces
are assumed carried out in the compactly generated category.
It is in this sense that the realization $|K|$
is always a topological group;
in general, the multiplication map will be continuous only
in the compactly generated refinement of the product topology
on $|K| \times |K|$.
For countable $K$, there is no difference, though.
Here is our main result.

\proclaim{Theorem}
There is a 
canonical
homeomorphism
of principal $|K|$-bundles
between 
the realization $|WK| \to |\overline WK|$ 
of the $W$-construction
and 
the classifying bundle
$E|K| \to B|K|$
which is natural in $K$.
\endproclaim

The map from 
$|WK|$ to $E|K|$
could be viewed as a kind of
perturbed
geometric Alexander-Whitney map 
while the map in the other direction
is a kind of
perturbed
geometric shuffle map 
(often referred to as Eilenberg-Zilber map)
but this analogy should not
be taken too far.
\smallskip
The classifying space
$B|K|$
is the realization
of the nerve $NK$ of $K$
as a {\it bisimplicial\/} set.
The latter is homeomorphic
to the realization of 
its {\it diagonal\/} $DNK$
since this is known to be true for an arbitrary
bisimplicial set \cite\quillone.
The diagonal $DNK$,
in turn, does {\it not\/} coincide
with the reduced $W$-construction
$\overline WK$, though, but after realization the two are homeomorphic.
We shall spell out the precise relationships in Section 4 below.
\smallskip
Eilenberg-Mac Lane
introduced the bar and $W$-constructions 
in \cite\eilmactw\ 
and
showed that,
for any (connected) simplicial algebra
$A$,
there is a \lq\lq reduction\rq\rq\ 
of the realization $|\overline WA|$ of
the reduced $W$-construction of $A$ onto
the (reduced normalized) bar construction
$B|A|$ 
of the normalized chain algebra $|A|$ of $A$
and raised the question whether
this reduction is in fact part of a {\it contraction\/}.
By means of homological perturbation theory,
in his \lq\lq Diplomarbeit\rq\rq\ \cite\wongone\ 
supervised by the second named author,
Wong
answered
this question
by establishing such a contraction.
Wong's basic tool 
is the \lq\lq perturbation lemma\rq\rq\ 
exploited in \cite\gugenhtw;
see \cite\huebkade\ for details and history.
\smallskip
Our result,
apart from being interesting in its own right,
provides a step towards a rigorous understanding
of lattice gauge theory.
See \cite\lattice\ for details.
Using the notation 
$K_Y$
for Kan group 
{\rm \cite\kanone}\
of a reduced simplicial set $Y$,
at this stage, we only spell out the following
consequence, relevant for
what is said in the references just cited.

\proclaim{Corollary}
For a reduced
simplicial set $Y$, there is a canonical map
from its realization $|Y|$ to
the classifying space
$B|K_Y|$
of the realization of 
$K_Y$
which is natural in $Y$ and a homotopy equivalence.
\endproclaim

The proof of our main result involves certain recursive descriptions
of the $W$-construction and classifying bundle
and relies on the facts that the realization functor
carries an action of a simplicial group to a geometric
action of its realization and preserves
reduced cones and colimits.
It would be interesting to extend the results of the present paper
to simplicial groupoids, so that a result of the kind 
given in the Corollary
would follow for an arbitrary connected simplicial set,
with the Kan group replaced by the Kan groupoid \cite\dwykanon.
We hope to return to this issue elsewhere.
The argument we shall give here involves suitable (reduced)
cone constructions
which no longer serve for the desired purpose  for the more
general groupoid case.
\smallskip
We are indebted to Jim Stasheff and to the referee for a number
of most helpful comments.

\beginsection 1. The classifying space of a topological group

Let $G$ be a topological group.
Its {\it nerve\/}
$
NG
$
\cite\bottone, \cite\botshust,
\cite\gsegatwo\ 
is the simplicial space
having in degree $k\geq 0$ the constituent $NG_k = G^{\times k}$, 
with the standard simplicial operations.
The usual {\it lean\/}
realization
$BG=|NG|$ of 
$NG$
is a classifying space for $G$,
cf. 
\cite\maclafou,
\cite\gsegatwo,
\cite\stasheig;
there is an analoguous
construction of contractible  total space
$EG$ together with a free $G$-action and projection
$\xi$ onto $BG$,
and this projection
is locally trivial provided
$(G,e)$ is an NDR 
(neighborhood deformation retract)
\cite\steenone.
We note, for completeness,
that the {\it fat\/}
realization
$||NG||$
yields {\smc Milnor's} 
classifying space
\cite\milnoru,
and the projection from the corresponding total space
to $||NG||$
is always locally trivial
whether or not
$(G,e)$ is an  NDR.
Below
$(G,e)$ will always be 
CW-pair and hence
an NDR, cf. e.~g. the discussion in the appendix to \cite\gsegathr,
and we shall exclusively deal with the lean realization
$BG =|NG|$.
To reproduce a description thereof, and to introduce notation,
write 
$\Delta$ for the category  of finite ordered sets
$[q] = (0,1,\dots,q),\, q\geq 0,$
and monotone maps.
We recall the standard {\it coface\/} and {\it codegeneracy\/} operators 
$$
\align
\varepsilon^j
&\colon [q-1] @>>> [q],
\quad
(0,1,\dots, j-1,j,\dots,q-1)
\mapsto (0,1,\dots, j-1,j+1,\dots,q),
\\
\eta^j
&\colon [q+1] @>>> [q],
\quad
(0,1,\dots, j-1,j,\dots,q+1)
\mapsto (0,1,\dots, j,j,\dots,q),
\endalign
$$
respectively.
As usual, 
for a simplicial object,
the corresponding face and degeneracy operators
will be written $d_j$ and $s_j$.
The assignment 
to
$[q]$ of the standard simplex
$\nabla [q]=\Delta_q$
yields a {\it cosimplicial\/}
space $\nabla$;
here we wish to distinguish
clearly
in notation between
the cosimplicial space
$\nabla$
and the category $\Delta$.
The lean  geometric realization
$|NG|$ is
the {\it coend\/}
$NG \otimes _\Delta\nabla$,
cf. e.~g. \cite\maclbotw\ 
for details on this notion.
Exploiting this observation,
Mac Lane 
observed
in \cite\maclafou\ 
that
$|NG|$
coincides with the classifying space
for $G$ constructed by
Stasheff \cite\stashone\ 
and Milgram \cite\milgrone;
see also
Section 1 of Stasheff's survey paper
\cite\stasheig\ 
and Segal's paper \cite\gsegatwo.
Mac Lane actually worked with a 
variant of the category
$\Delta$ which enabled him to
handle simultaneously
the total space $EG$ and the base $BG$.
\smallskip
Steenrod \cite\steenone\ 
has given a recursive description
of
$|NG|$
which we shall subsequently use.
For ease of exposition,
following [1],
we reproduce it briefly
in somewhat more categorical language.
This will occupy the rest of this Section.
\smallskip
For a space
$X$ endowed with a $G$-action $\phi \colon X \times G \to X$
we write
$
\eta = \eta^G_X
\colon 
X \to X \times G$
for
the {\it unit\/}
given by
$\eta(x) = (x,e)$.
For an arbitrary space $Y$,
right translation of $G$ induces an obvious free $G$-action
$\mu$ on $Y \times G$.
In 
categorical language
\cite\maclbotw,
the 
functor $\times G$ and natural transformations $\mu$ and $\eta$
constitute a {\it monad\/} $(\times G,\mu,\eta)$
 and a $G$-action
on a space $X$ is an {\it algebra\/} structure on $X$
over this monad.
Sometimes we shall refer to an action of a topological group on a space
as a {\it geometric\/} action.
\smallskip
Let $D$ be any space and $E$ a subspace
endowed with a $G$-action
$\phi \colon E \times G \to E$;
the inclusion of $E$ into $D$ is written $\beta$.
Recall that the {\it enlargement\/}
$\overline D \supseteq D$
of the $G$-action
is characterized by the property:
if $Y$ is any $G$-space, and $f$ any map from
$D$ to $Y$ whose restriction to $E$ is a $G$-mapping, then there
exists a unique $G$-mapping
$\overline f$ from $\overline D$ to $Y$ extending $f$.
The space 
$\overline D$
then fits into a push out diagram
$$
\CD
E \times G
@>\phi>>
E
\\
@V{\beta \times \roman{Id}}VV
@VVV
\\
D \times G
@>>>
\overline D
\endCD
\tag1.1
$$
and this provides a construction for $\overline D$.
Moreover, right action of $G$ on $D \times G$
induces an action 
$$
\overline  \phi \colon \overline D \times G @>>> \overline D
\tag1.2
$$
of $G$ on $\overline D$, and the composite 
$$
\alpha \colon D @>>>
\overline D
\tag1.3
$$
of the unit
$\eta \colon D \to D \times G$
with the map from
$D\times G$
to
$\overline D$
in (1.1)
embeds
$D$
into $\overline D$.
When $D$ is based and $E$ is a based subspace, the products 
$E \times G$
and
$D \times G$
inherit an obvious base point,
and the square (1.1) is one in the category of based spaces
whence, in particular,
the enlargement
$\overline D$ inherits a base point.
This notion of enlargement
of $G$-action is functorial in the appropriate sense.
See \cite\steenone\ for details.
This kind of universal construction is available whenever
one is given an algebra structure over a monad preserving
push out diagrams.
\smallskip
The unit interval $I = [0,1]$
is a topological monoid under ordinary multiplication
having $1$ as its unit,
and hence we can talk about
an $I$-action $X \times I \to X$
on a space $X$.
Such an $I$-action
is plainly a special kind of homotopy
which, for $t=1$, is the identity.
In the above categorical spirit,
the interval $I$ gives rise to a monad
$(\times I,\mu,\eta)$
and an $I$-action
on a space $X$ is an {\it algebra\/} structure on $X$
over this monad.
\smallskip
The {\it base point\/} of $I$ is defined to be 0.
Following \cite\steenone,
for a based space $(X,x_0)$,
we shall 
refer to an $I$-action
$\psi \colon X \times I \to X$
as a {\it contraction\/}
of $X$ 
({\it to the base point\/} $x_0\in X$)
provided $\psi$ 
sends the base point $(x_0,0)$ of $X \times I$ to $x_0$ and
factors through
the {\it reduced cone\/}
or {\it smash product\/}
$$
CX = X\wedge I = X \times I \big/( X \times \{0\} \cup \{x_0\} \times I)
$$
that is to say,
$$
\psi(x,0) = x_0 = \psi(x_0,t)
$$
for all $x \in X,\, t \in I$;
the reduced cone
will be endowed with the obvious base point,
the image of $X \times \{0\} \cup \{x_0\} \times I$
in $CX$.
Whenever we say \lq\lq contraction\rq\rq, we mean \lq\lq contraction to
a pre-assigned base point\rq\rq.
Abusing notation,
the corresponding map from
$CX$
to $X$ will as well be denoted by $\psi$
and referred to as a {\it contraction\/}.
Moreover we write
$\eta = \eta^C_X$
for the map,
the corresponding {\it unit\/},
which embeds $X$ into $CX$ by sending
a point $x$ of $X$ to $(x,1) \in CX$.
The right action of $I$ on $X \times I$
induces a contraction  $\mu^C_X\colon CCX \to CX$
of $CX$.
Again we can express this in 
categorical language:
the 
functor $C$ and natural transformations
$\mu$ and $\eta$
constitute a {\it monad\/} 
and a contraction
of a based space $X$ is an {\it algebra\/} structure on $X$
over this monad.
Sometimes we shall refer to a contraction of a space
as a {\it geometric\/} contraction.
\smallskip
Let $(E,x_0)$ be any based space and $(D,x_0)$ a based subspace
endowed with a contraction
$\psi\colon CD\to D$;
the inclusion of $D$ into $E$ is written $\alpha$.
The {\it enlargement\/}
$(\overline E, x_0) \supseteq (E,x_0)$
of the contraction
is characterized by the property:
if  $f$ is any map from
$E$ 
to a space $Y$ having a contraction to some point $y_0$
whose restriction to $D$ is an $I$-mapping, then there
exists a unique $I$-mapping
$\overline f$ from $\overline E$ to $Y$ extending $f$.
The space  
$\overline E$
then fits into a push out diagram
$$
\CD
CD
@>\psi>>
D
\\
@V{C\alpha}VV
@VVV
\\
CE
@>>>
\overline E
\endCD
\tag1.4
$$
which provides a construction for
$\overline E$.
Moreover,
the composite 
$$
\beta \colon E @>>>
\overline E
\tag1.5
$$
of the unit
$\eta \colon E \to CE$
with the map from
$CE$
to
$\overline E$
in (1.4)
embeds
$E$
into $\overline E$
and
the right action of $I$ on $E \times I$
induces a contraction of $CE$ which, in turn,
induces a contraction 
$$
\overline \psi \colon C\overline E @>>> \overline E
\tag1.6
$$
of $\overline E$.
This notion of enlargement
of contraction is functorial in the appropriate sense.
See \cite\steenone\ for details.
\smallskip
Alternating the above constructions,
in \cite\steenone,
Steenrod defines based spaces and injections
of based spaces
$$
D_0 @>{\alpha_0}>> E_0 @>{\beta_0}>> D_1 @>{\alpha_1}>> 
\dots @>{\beta_{n-1}}>> D_n @>{\alpha_n}>> E_n
@>{\beta_n}>> D_{n+1} @>{\alpha_{n+1}}>> \dots
\tag1.7
$$
by induction on $n$
together with
contractions
$\psi_n\colon CD_n \to D_n$
(Steenrod writes these contractions as $I$-actions
$D_n \times I \to D_n$)
and $G$-actions
$\phi_n \colon E_n \times G \to E_n$
in the following way:
Let $D_0$ 
consist
of the single point $e$ with the obvious contraction.
Let $E_0 = G$, the right action being right translation.
Now define
$(D_1,e)$ to be the enlargement
to $(E_0,e)$,
$(\overline E_0,e)$, of the contraction of
$(D_0,e)$;
then $D_1$ is just the reduced cone on $E_0$.
Define $E_1$ to be the enlargement
to $D_1$,
$\overline D_1$,
of the $G$-action on $E_0$.
In general, $D_n$ is the enlargement to
$(E_{n-1},e)$,
$(\overline E_{n-1},e)$,
 of the contraction
$\psi_{n-1}$ of $(D_{n-1},e)$
so that $D_n$ fits into a push out
square
$$
\CD
CD_{n-1}
@>{\psi_{n-1}}>>
D_{n-1}
\\
@V{C\alpha_{n-1}}VV
@VVV
\\
CE_{n-1}
@>>>
D_n;
\endCD
\tag1.8
$$
the requisite injection $\beta_{n-1} \colon E_{n-1} \to D_n$
is the map denoted above by $\beta$,
cf. (1.5);
and
the requisite contraction
$\psi_{n}\colon CD_n @>>> D_n$
of $(D_n,e)$
or, equivalently, $I$-action
$\psi_{n}\colon D_n  \times I @>>> D_n$,
is the map denoted above by $\overline\psi$,
cf. (1.6).
Likewise,
$E_n$ is the enlargement
to $D_n$,
$\overline D_n$,
of the $G$-action $\phi_{n-1}$ on $E_{n-1}$,
so that $E_n$ fits into a push out square
$$
\CD
E_{n-1} \times G
@>{\phi_{n-1}}>>
E_{n-1}
\\
@V{\beta_{n-1} \times \roman{Id}}VV
@VVV
\\
D_n \times G
@>>>
E_n;
\endCD
\tag1.9
$$
the requisite
$G$-action
$\phi_n \colon E_n \times G
@>>>
E_n
$
and
injection $\alpha_n \colon D_n \to E_n$
are the action denoted above by $\overline \phi$,
cf. (1.2),
and the map denoted above by $\alpha$,
cf. (1.3),
respectively.
The union
$$
E_G = \bigcup_{n=0}^{\infty} E_n = \bigcup_{n=0}^{\infty} D_n,
$$
endowed with the weak topology,
inherits
a $G$-action
$\phi \colon E_G \times G \to E_G$
and contraction
$\psi \colon CE_G \to E_G$
from
the $\phi_n$'s 
and
$\psi_n$'s,
respectively.
The $G$-action is free,
and the orbit space
$BG = E_G \big / G$
equals the lean geometric realization
$|NG|$ of the nerve of $G$.
This is Steenrod's result in \cite\steenone.

\beginsection 2. The recursive description of the $W$-construction

In \cite\cbergthe,
the first named author observed
that 
the $W$-construction
admits 
a recursive description
of formally the same kind
as (1.7) above,
except
that it is carried out in the category of based simplicial sets.
This is among the key points of the paper.
We shall explain it now.
\smallskip
Let $K$ be a simplicial group.
Let $e$ denote the trivial simplicial group
viewed
at the same time
as the simplicial point.
For a simplicial set
$X$ endowed with a $K$-action $\phi \colon X \times K \to X$
we write
$
\eta = \eta^K_X
\colon 
X \to X \times K$
for
the {\it unit\/}
of the action;
in each degree, it is given by
$\eta(x) = (x,e)$.
Given an arbitrary simplicial set $Y$,
right translation of $K$ induces an obvious action
$\mu$ of $K$ on $Y \times K$.
Much as before,
in 
categorical language,
the 
functor $\times K$ and natural transformations $\mu$ and $\eta$
constitute a {\it monad\/} $(\times K,\mu,\eta)$
in the category of simplicial sets
and a $K$-action
on a simplicial set $X$ is an {\it algebra\/} structure on $X$
over this monad.
Moreover realization preserves
monad and algebra structures.
In other words:
the realization
of 
a $K$-action
$\phi \colon X \times K \to X$
on a simplicial set $X$
is a
geometric action
$|\phi| \colon |X| \times |K| \to |X|$
in the usual sense.
Notice this 
involves the standard homeomorphism 
\cite\milnothr\ 
between the
realization 
$|X \times K|$
of
the simplicial set
$X \times K$
and the product 
$|X| \times |K|$
of the
realizations
(with the compactly generated topology).
The homeomorphism 
between
$|X \times K|$
and
$|X| \times |K|$
is of course natural and relies
on the fact that,
for an arbitrary
bisimplicial set,
the realization of the diagonal is homeomorphic
to the realization as a bisimplicial set,
cf. \cite\quillone\ (Lemma on p. 86).
Note, however, that
the product
$|X| \times |K|$
yields a realization of
$X \times K$
only after subdivision of
the product CW-decomposition of
$|X| \times |K|$, cf. \cite\puppeone\ (Satz 5 p. 388).
This reflects, of course, the decomposition
coming into play in the Eilenberg-Zilber map.
\smallskip
Recall that in the category of simplicial sets there are {\it two\/}
natural (reduced) cone constructions. The first one is defined by the
simplicial smash product with the standard simplicial model
$\Delta[1]$ of the unit interval.
We shall say more about this in Section 4 below.
The recursive description of the $W$-construction
crucially involves the 
second 
somewhat more economical
cone construction
which relies on the observation that 
an $(n+1)$-simplex serves as
a cone on an $n$-simplex.
We reproduce this cone construction briefly;
it differs from the
one given in \cite\curtione\ (p. 113)
by the order of face and degeneracy operators;
our convention is forced
here by our description of
the $W$-construction with structure group acting from the right,
cf. what is said in (2.6) below.
\smallskip
Let
$X$ be a simplicial set.
For $j \geq 0$,
we shall need countably many disjoint copies
of each $X_j$ which
we describe in the following way:
For $j \geq 0$, 
consider the cartesian product
$X_j \times \bold N$ with the natural numbers $\bold N$.
Let $o$ be a point which we formally 
assign dimension -1 and,
given $i \in \bold N$, 
write
$X_{-1}(i) = 
\{(o,i)\}$
so that
each $X_{-1}(i)$
consists of a single element;
next, for $j \geq 0$, let
$X_j(i) = X_j \times \{i\}$.
The 
{\it unreduced simplicial
cone\/}
$\widehat CX$
on $X$ is given by
$$
(\widehat CX)_n = X_n(0) \cup \dots \cup X_0(n) \cup X_{-1}(n+1),
\quad
n \geq0,
$$
with face and degeneracy operators given by the formulas
$$
\aligned
d_j(x,i) &= \cases (d_jx,i) \quad & j \leq n-i \\
                   (x,i-1)  \quad & j > n-i 
            \endcases
\\
s_j(x,i) &= \cases (s_jx,i) \quad & j \leq n-i \\
                   (x,i+1)  \quad & j > n-i 
            \endcases
\endaligned
$$
Notice that 
in these formulas $n-i= \dim x$; in particular,
$$
d_j(o,n+1) = (o,n),
\qquad 
s_j(o,n) = (o,n+1),
\qquad 0 \leq j \leq n.
$$
\smallskip
Let now
$(X,*)$
be a {\it based\/} simplicial set.
The unreduced simplicial cone
$\widehat C\{*\}$
of the simplicial point 
$\{*\}$
is the simplicial interval,
and
the {\it reduced simplicial
cone\/} $CX$
is simply the quotient
$$
CX = \widehat CX \big / \widehat C\{*\}.
$$
For each $n \geq 0$,
its constituent $(CX)_n$ arises
from
the union 
$X_n(0) \cup \dots \cup X_0(n)$
by identifying
all $(*,i)$ to a single point written $*$,
the {\it base point\/} of $CX$.
The non-degenerate simplices of $CX$ different from the base point
look like $(x,0)$
and $(x,1)$
where $x$ runs through non-degenerate simplices of $X$.
We
write $\eta = \eta^C_X \colon X \to CX$
for the {\it unit\/}
induced by the assignment to
$x \in X_n$ of $(x,0) \in X_n(0)$.
A (simplicial) contraction is, then, a
morphism
$\psi \colon CX \to X$
of based simplicial sets
satisfying
$$
\psi \circ \eta = \roman{Id}_X.
$$
The cone $CX$ itself admits the obvious contraction
$$
\mu =\mu^C_X\colon CCX \to CX,
\quad
((x,i),j) \mapsto (x,i+j).
$$
A contraction $\psi$ is called {\it conical\/}
provided
$$
\psi \circ C\psi = \psi \circ \mu.
$$
We note that 
a geometric contraction in the sense of
Section 1 above,
 being defined as an action of
the associative monoid $I$,
automatically
satisfies the usual
associativity 
law 
for an action.
Under the present circumstances,
the property of being conical corresponds
to this associativity property.
The contraction $\mu^C_X$ of $CX$ is conical;
in categorical terms, 
the triple $(C,\mu,\eta)$ is a monad
in the category of simplicial sets, and
a conical contraction
is an algebra structure 
in the category of simplicial sets
over this monad.
\smallskip
We now have the machinery in place
to reproduce
the first named author's
recursive description
the $W$-construction:
Define
based simplicial sets and injections of based
simplicial sets
$$
D_0 @>{\alpha_0}>> E_0 @>{\beta_0}>> D_1 @>{\alpha_1}>> 
\dots @>{\beta_{n-1}}>> D_n @>{\alpha_n}>> E_n
@>{\beta_n}>> D_{n+1} @>{\alpha_{n+1}}>> \dots
\tag2.1
$$
by induction on $n$
together with
conical contractions
$\psi_n\colon CD_n \to D_n$
and $K$-actions
$\phi_n \colon E_n \times K \to E_n$
 on each $E_n$
from the {\it right\/}
in the following way:
Let $D_0 = e$, with the obvious conical contraction $\psi_0$,
let $E_0 = K$, 
viewed as a based simplicial set in the obvious way,
the right action $\phi_0$ being translation,
and let
$\alpha_0$ be the obvious morphism of based simplicial sets
from
$D_0$ to $E_0$.
For $n \geq 1$, define
$(D_n,e)$ to be the {\it enlargement\/} to
$(E_{n-1},e)$ of the contraction
$\psi_{n-1}\colon CD_{n-1} \to D_{n-1}$, 
that is,
$D_n$ is characterized by the requirement that
the diagram
$$
\CD
CD_{n-1}
@>{\psi_{n-1}}>>
D_{n-1}
\\
@V{C\alpha_{n-1}}VV
@VVV
\\
CE_{n-1}
@>>>
D_n
\endCD
\tag2.2
$$
be a push out square of (based) simplicial sets;
the 
composite of the unit 
$\eta$ from
$E_{n-1}$
to
$CE_{n-1}$
with the morphism
$CE_{n-1}
@>>>
D_n$
of simplicial sets
in (2.2)
yields the
requisite injection $\beta_{n-1} \colon E_{n-1} \to D_n$,
and
the  contraction
$\psi_{n-1}$ and the conical
contraction of $CE_{n-1}$ induce a conical contraction
$\psi_{n}\colon CD_n @>>> D_n$.
Likewise,
$E_n$ is the {\it enlargement\/}
to $D_n$
of the $K$-action $\phi_{n-1}$ on $E_{n-1}$,
that is,
$E_n$ is characterized by a push out square
of based simplicial sets of the kind
$$
\CD
E_{n-1} \times K
@>{\phi_{n-1}}>>
E_{n-1}
\\
@V{\beta_{n-1} \times \roman{Id}}VV
@VVV
\\
D_n \times K
@>>>
E_n ;
\endCD
\tag2.3
$$
the requisite
$K$-action
$\phi_n \colon E_n \times K
@>>>
E_n
$
is 
induced by $\phi_{n-1}$ and
the obvious 
$K$-action on
$D_n \times K$, and
the requisite injection $\alpha_n \colon D_n \to E_n$
is the 
composite of the unit with the morphism
$D_n \times K
@>>>
E_n$
of simplicial sets
in (2.3).
The limit
$$
WK = \lim_{\to} E_n = \lim_{\to} D_n
$$
inherits
a $K$-action
$\phi \colon WK \times K \to WK$
and conical contraction
$\psi \colon CWK \to WK$
from
the $\phi_n$'s 
and
$\psi_n$'s,
respectively.
The $K$-action is free,
and the projection map to the quotient
$\overline WK = WK \big / K$
yields the universal simplicial $K$-bundle
$$
WK
@>>>
\overline WK
$$
or $W$-construction of $K$,
cf. \cite\cbergthe,
with action of $K$ from the {\it right\/}.
\smallskip
For intelligibility,
we 
explain some of the requisite details:
A straightforward induction establishes the following descriptions
of the simplicial sets $D_k$ and $E_k$:
$$
\align
(D_k)_n &=
\left \{
(i_0, k_0,i_1,k_1,\dots, k_{\ell -1}, i_\ell)\ | \ 
0 \leq \ell\leq k,\ i_s \geq 0,
\right.
\\
&\ \left.
n = i_0 + \dots + i_\ell,\ k_s \in K_{i_0 + \dots + i_s},
\ 0 \leq s <\ell \right\}
\big/ \sim ,
\\
(E_k)_n &=
\left \{
(i_0, k_0,i_1,k_1,\dots, k_{\ell -1}, i_\ell)\ | \ 
0 \leq \ell\leq k,\ i_s \geq 0,
\right.
\\
&\ \left.
n = i_0 + \dots + i_\ell,\ k_s \in K_{i_0 + \dots + i_s},
\ 0 \leq s \leq \ell \right\}
\big/ \sim ,
\endalign
$$
where
$$
(\dots,i_s,e,i_{s+1},\dots)
\sim
(\dots,i_s + i_{s+1},\dots),
\quad
(\dots,k_s,0,k_{s+1},\dots)
\sim
(\dots,k_s k_{s+1},\dots).
$$
Thus,
for $n \geq 0$,
$$
\align
(WK)_n &=
\left \{
(k_{j_0}k_{j_1} 
\dots
k_{j_\ell}\ | \ 
0 \leq j_0 < \dots < j_\ell =n\quad\text{and}\quad
\right.
\\
&\ \left.
k_{j_s} \in K_{j_s} \setminus  e_{j_s},\ 0 \leq s <\ell,
\ k_{j_\ell} \in K_{j_\ell}
 \right\}.
\endalign
$$
{}From this,
adding the requisite neutral elements
wherever appropriate,
we deduce
the following 
more common
explicit description:
For $n \geq 0$,
$$
(WK)_n = K_0 \times \dots \times K_n,
$$
with face and degeneracy operators given by the formulas
$$
\aligned
d_0(x_0,\dots, x_n) 
&= (d_0x_1,\dots, d_0 x_n)\\
d_j(x_0,\dots, x_n) 
&=(x_0,\dots,x_{j-2}, x_{j-1}d_j x_j, d_j x_{j+1},\dots,d_j x_n),
\quad 1 \leq j \leq n
\\
s_j(x_0,\dots, x_n) 
&=(x_0,\dots,x_{j-1}, e,s_j x_j,s_j x_{j+1},\dots,s_j x_n),
\quad 0 \leq j \leq n;
\endaligned
\tag2.4
$$
further,
$(\overline WK)_0 = \{e\}$ and, for $ n \geq 1$
$$
(\overline WK)_n = K_0 \times \dots \times K_{n-1},
$$
with face and degeneracy operators given by the formulas
$$
\aligned
d_0(x_0,\dots, x_{n-1}) 
&= (d_0x_1,\dots, d_0 x_{n-1}),
\\
d_j(x_0,\dots, x_{n-1}) 
&=(x_0,\dots,x_{j-2}, x_{j-1}d_j x_j, d_j x_{j+1},\dots,d_j x_{n-1}),
\\
&\qquad 1 \leq j \leq n-1,
\\
d_n(x_0,\dots, x_{n-1}) 
&=(x_0,\dots, x_{n-2}), 
\\
s_0(e) &= e 
\in K_0,
\\
s_j(x_0,\dots, x_{n-1}) 
&=(x_0,\dots,x_{j-1}, e,s_j x_j,s_j x_{j+1},\dots,s_j x_{n-1}),
\\
&\qquad 0 \leq j \leq n.
\endaligned
\tag2.5
$$
\smallskip\noindent
{\smc Remark 2.6.}
Here preferred treatment is given to the {\it last\/} face operator, as is
done in 
\cite\gugenmay\ and
\cite\kanone.
This turns out to be the appropriate thing to do
for principal bundles with  structure group acting
on the  total space
from the {\it right\/}
and simplifies comparison with the bar construction.
See for example what is said on p. 75
of \cite\gugenmay.
The formulas (2.4) and (2.5) arise from those given in (A.14) of
\cite\gugenmay\ for a simplicial algebra
by the obvious translation to the corresponding formulas
for a simplicial monoid;
they differ from those in \cite\curtione\ (pp. 136 and 161)
where the constructions are carried out
with structure group acting from the {\it left\/}.

\beginsection 3. The proof of the Theorem

The realization of a conical contraction
$\psi \colon CX \to X$
of a based simplicial set $(X,x_0)$
is
a geometric contraction
$|\psi| \colon C|X| \to |X|$
in the sense reproduced in Section 1 above.
In fact,
the 
association
$$
(|x|(t_0,\dots,t_n),t)
\longmapsto
(|(x,1)|(tt_0,\dots,tt_n,1-t),
\ x \in X_n,\ n \geq 0,
$$
yields a homeomorphism from
the reduced cone $C|X|$
on the realization $|X|$
to the realization $|CX|$
of the cone
and, furthermore,
the realizations of the unit $\eta$
and
$C$-algebra structure
$\mu^C_X\colon CCX \to CX$
yield the geometric unit
$|X| \to C|X|$
and
geometric $C$-algebra structure
$\mu^C_{|X|}\colon CC|X| \to C|X|$,
that is, the realization
preserves monad- and $C$-algebra structures.
\smallskip
The proof of the Theorem is now
merely an elaboration of the observation that
the realization functor
$|\cdot|$
carries an action 
of a simplicial group to a geometric action of its realization,
preserves reduced cones
and,
having a right adjoint (the singular complex functor),
also preserves  colimits.
In fact,
denote the corresponding sequence (1.7) 
of based topological spaces
for the realization $|K|$ by
$$
D_0|K| @>{\alpha_0 |K|}>> E_0|K| @>{\beta_0|K|}>> 
\dots
@>{\alpha_n |K|}>> E_n|K|
@>{\beta_n |K|}>> D_{n+1}|K| 
@>{\alpha_{n+1}|K|}>> 
\dots
\tag3.1
$$
and, likewise,
write
$$
D_0 K  @>{\alpha_0 K}>> E_0 K  @>{\beta_0 K}>> 
\dots
@>{\alpha_n K}>> E_n K 
@>{\beta_n K}>> D_{n+1} K  
@>{\alpha_{n+1} K}>> 
\dots
\tag3.2
$$
for the corresponding sequence (2.1)
in the category of based simplicial sets.
Realization carries
the sequence (3.2) to the sequence 
$$
|D_0 K|  @>{|\alpha_0 K|}>> |E_0 K|  @>{|\beta_0 K|}>> 
\dots
@>{|\alpha_n K|}>> |E_n K| 
@>{|\beta_n K|}>> |D_{n+1} K|  
@>{|\alpha_{n+1} K|}>> 
\dots
\tag3.3
$$
of based topological spaces.
Now
$$
D_0 |K|= e=|D_0 K|,\ E_0 |K| = |K|  =|E_0 K|,
$$
and the map $\alpha_0 |K| = |\alpha_0 K|$
is the canonical inclusion. 
Let 
$$
\tau_0 \colon D_0 |K| @>>> |D_0 K|
\quad \text{and}\quad                   
\rho_0 \colon E_0 |K| @>>> |E_0 K|
$$
be the identity mappings.
Let $n \geq 1$ and
suppose by induction that homeomorphisms
$$
\tau_j \colon D_j |K| @>>> |D_j K|\quad
\text{and}\quad                   
\rho_j \colon E_j |K| @>>> |E_j K|,
$$
each $\rho_j$ being $|K|$-equivariant,
have been constructed for
$j < n$,
having the following properties:
\roster
\item
The diagrams
$$
\CD
E_{j-1}|K|
@>{\beta_{j-1}|K|}>>
D_j|K|
@.
{\qquad \qquad}
@.
D_j|K|
@>{\alpha_j|K|}>>
E_j|K|
\\
@V{\rho_{j-1}}VV
@VV{\tau_j}V
@.
@V{\tau_j}VV
@VV{\rho_j}V
\\
|E_{j-1} K|
@>{|\beta_{j-1} K|}>>
|D_jK|
@.
{\qquad \qquad}
@.
|D_j K|
@>{| \alpha_j K|}>>
|E_j K|
\endCD
\tag3.4
$$
are commutative;
\item
each $\tau_j$ identifies the 
realization
$|\psi_j K|\colon |CD_j K| @>>> |D_j K|$
of the conical
contraction
$\psi_j K\colon CD_jK @>>> D_j K$
of simplicial sets with the geometric contraction
$\psi_j |K|\colon CD_j |K| @>>> D_j |K|$;
\item
each
$\rho_j$ 
identifies the 
realization
$|\phi_j K|\colon |(E_j K) \times K | @>>> |E_j K|$
of the 
simplicial $K$-action
$\phi_j K\colon (E_jK) \times K @>>> E_j K$
with the topological $|K|$-action
$\phi_j |K|\colon (E_j |K|) \times |K| @>>> E_j |K|$.
\endroster
Consider the realization of (2.2);
it is a push out square of topological spaces.
Hence 
the 
maps $\tau_{n-1}$ and $\rho_{n-1}$
induce a map
$\tau_n$ from
$D_n|K|$ to
$|D_n K|$, necessarily a homeomorphism,
so that
$C|\tau_{n-1}|,|\tau_{n-1}|, C|\rho_{n-1}|$ and
$|\tau_n|$
yield a homeomorphism
of squares
between the realization
of
(2.2) and (1.8).
Moreover,
the homeomorphism
$\tau_n$ identifies the 
realization
$|\psi_n K|\colon |CD_n K| @>>> |D_n K|$
of the conical
contraction
$\psi_n K\colon CD_nK @>>> D_n K$
of simplicial sets with the contraction
\linebreak
$\psi_n |K|\colon CD_n |K| @>>> D_n |K|$.
Likewise
the 
maps $\rho_{n-1}$ 
and $\tau_{n}$
induce a map
$\rho_n$ from
$E_n|K|$ to
$|E_n|K|$, necessarily a $|K|$-equivariant homeomorphism,
so that
$|\rho_{n-1}\times \roman{Id}_K |,|\rho_{n-1}|, |\tau_n\times \roman{Id}_K |$ 
and
$|\rho_n|$
yield a homeomorphism
of squares
between the realization
of
(2.3) and (1.9).
Moreover,
the homeomorphism
$\rho_n$ 
is $|K|$-equivariant and
identifies the 
realization
$|\phi_n K|\colon |(E_n K) \times K | @>>> |E_n K|$
of the simplicial $K$-action
$\phi_n K\colon (E_nK) \times K @>>> E_n K$
with the topological $|K|$-action
$\phi_n |K|\colon (E_n |K|) \times |K| @>>> E_n |K|$.
The requisite diagrams
(3.4) 
for $j=n$
are manifestly commutative.
This completes the inductive step.
\smallskip
The limit
$$
\rho = \lim \rho_n = \lim \tau_n
\colon E|K| @>>> |WK|
$$
is a $|K|$-equivariant homeomorphism;
it 
identifies the principal
$|K|$-bundles
\linebreak
$E|K| \to B|K|$ and
$|WK| \to |\overline WK|$ 
as asserted, is plainly
natural in $K$
and, in particular,
induces a natural homeomorphism
from
$B|K|$ to $|\overline WK|$.
This proves the Theorem.

\beginsection 4. The other cone construction

The classifying space
$B|K|$
is the realization
of the nerve $NK$ of $K$
as a {\it bisimplicial\/} set.
On the other hand,
the diagonal $DNK$
is a simplicial set
which does {\it not\/} coincide
with the reduced $W$-construction
$\overline WK$, but its
realization 
is homeomorphic
to the realization of 
the nerve $NK$ of $K$
as a bisimplicial set
since this is known to be true for an arbitrary
bisimplicial set \cite\quillone.
The purpose of this Section is to clarify
the relationships
between the various spaces and constructions.
\smallskip
As already pointed out,
the construction (2.1) can be carried out with the simplicial smash product
$(\cdot) \wedge \Delta[1]$ instead of the reduced
cone:
The simplicial interval 
$\Delta[1]$
carries a (unique) structure of a simplicial monoid
having $(1)$ as its unit,
and hence we can talk about
an action $X \times \Delta[1] \to X$
of $\Delta[1]$
on a simplicial set $X$;
such an action
is a special kind of simplicial homotopy
which \lq\lq ends\rq\rq\ at the identity morphism of $X$.
The fact that the naive notion of homotopy
of 
morphisms of
simplicial sets is not an equivalence relation
is not of significance here.
Much as before, the 
simplicial interval 
$\Delta[1]$
 gives rise to a monad
$(\times \Delta[1],\mu,\eta)$
in the category of simplicial sets
and an action
of
$\Delta[1]$
on a simplicial set $X$ is an {\it algebra\/} structure on $X$
over this monad.
\smallskip
The {\it base point\/} of $\Delta[1]$ is defined to be $(0)$.
For a based simplicial set $(X,x_0)$,
we shall 
refer to an action
$\psi \colon X \times \Delta[1] \to X$
as a $\Delta[1]$-{\it contraction\/}
of $X$ 
provided $\psi$ 
sends the base point $(x_0,0)$ of $X \times \Delta[1]$ to $x_0$ and
factors through
the 
{\it simplicial  smash product\/}
$$
 X\wedge \Delta[1] 
= X \times \Delta[1] \big/( X \times \{0\} \cup \{x_0\} \times \Delta[1]).
$$
The latter
is viewed endowed with the obvious base point,
the image of $X \times \{0\} \cup \{x_0\} \times \Delta[1]$
in $ X\wedge \Delta[1]$.
Abusing notation,
the corresponding map from
$X\wedge \Delta[1]$
to $X$ will as well be denoted by $\psi$
and referred to as a $\Delta[1]$-{\it contraction\/}.
Moreover we write
$\eta = \eta^{\Delta[1]}_X$
for the map,
the corresponding {\it unit\/},
which embeds $X$ into $X\wedge \Delta[1]$ by sending
a simplex  $x$ of $X$ to $(x,1) \in X\wedge \Delta[1]$.
The right action of $\Delta[1]$ on $X \times \Delta[1]$
induces a $\Delta[1]$-contraction  
$$
\mu^{\Delta[1]}_X\colon X\wedge \Delta[1]\wedge \Delta[1] 
@>>> X\wedge \Delta[1]
$$
of $X\wedge \Delta[1]$.
In categorical language,
the 
functor $(\cdot) \wedge \Delta[1]$ and natural transformations
$\mu$ and $\eta$
constitute a {\it monad\/} 
in the category of simplicial sets,
and a $\Delta[1]$-contraction
of a based simplicial set $X$ is an {\it algebra\/} structure on $X$
over this monad.
\smallskip
Formally carrying out
the construction (2.1) with the simplicial smash product
$(\cdot) \wedge \Delta[1]$ instead of the reduced
cone
yields
based simplicial sets and injections of based
simplicial sets
$$
D'_0 @>{\alpha'_0}>> E'_0 @>{\beta'_0}>> D'_1 @>{\alpha'_1}>> 
\dots @>{\beta'_{n-1}}>> D'_n @>{\alpha'_n}>> E'_n
@>{\beta'_n}>> D'_{n+1} @>{\alpha'_{n+1}}>> \dots
\tag4.1
$$
together with
morphisms
$\psi'_n\colon D'_n\wedge \Delta[1] \to D'_n$
of simplicial sets having certain properties
and free $K$-actions
$\phi'_n \colon E'_n \times K \to E'_n$.
Its limit
$$
D = \lim_{\to} E'_n = \lim_{\to} D'_n
$$
inherits a morphism
$\psi'\colon D\wedge \Delta[1] \to D$
of simplicial sets
and a free $K$-action
$\phi'\colon D \times K \to D$.
To explain the significance thereof,
recall that the nerve construction yields a simplicial object
$$
K @>>> ENK @>>> NK
\tag4.2
$$
in the category of principal simplicial $K$-bundles
which is natural for morphisms of simplicial groups.
Here $ENK$ and
$NK$ inherit structures of bisimplicial sets,
one from the nerve construction and the other one from the simplicial
structure of $K$, and the projection from $ENK$ to $NK$ is a morphism
of bisimplicial sets;
further,
for each simplicial degree $q \geq 0$ coming from the
nerve construction,
(4.2) amounts to a principal $K$-bundle
$$
K_* @>>> (ENK)_{*,q} @>>> (NK)_{*,q}
$$
while 
for each simplicial degree $p \geq 0$
of $K$ = $\{K_p\}$ itself,
(4.2) comes down to the universal simplicial principal $K_p$-bundle
$$
K_p @>>> (ENK)_{p,*} @>>> (NK)_{p,*};
$$
in particular, each
$(ENK)_{p,*}$
is contractible in the usual sense.
The diagonal bundle
$$
\delta
\colon
DENK
@>>>
DNK
$$
is manifestly a principal $K$-bundle 
having $DENK$ contractible, and we have
$$
DENK = \lim_{\to} E'_n = \lim_{\to} D'_n
$$
as (right) $K$-set;
moreover, the above morphism
$\psi'\colon DENK \wedge \Delta[1] \to DENK$
induces a simplicial contraction of $DENK$.

\proclaim{Theorem 4.3}
There is a 
canonical
homeomorphism
of principal $|K|$-bundles
between 
the realization 
$|DENK| \to |\overline DNK|$ 
of the diagonal bundle
and 
the realization $|WK| \to |\overline WK|$ 
of the $W$-construction
which is natural in $K$.
\endproclaim

\demo{Proof}
The classifying space $B|K|$ is the realization
of $NK$ as a bisimplicial set, and the same kind of remark
applies to $E|K|$ and the projection to $B|K|$.
The already cited  fact that,
for an arbitrary
bisimplicial set,
the realization of the diagonal is homeomorphic
to the realization as a bisimplicial set \cite\quillone\ 
implies 
the following statement.

\proclaim{4.4}
There is a canonical $|K|$-equivariant homeomorphism
between $|DENK|$ and $E|K|$ and hence a canonical 
homeomorphism
between $|DNK|$ and $B|K|$.
These homeomorphisms are natural in $K$.
\endproclaim

We conclude from this that
{\it the statement of the Theorem\/}
(in the Introduction)
{\it is formally equivalent to the statement
of\/}
(4.3).
In fact, the Theorem identifies the realization 
the of $W$-construction with the
realization of the nerve as a {\it bisimplicial set\/}
whereas (4.3)
identifies the realization of the $W$-construction with the
realization of the {\it diagonal\/} of the nerve.
\enddemo

\smallskip\noindent
{\smc Remark 1.}
While the statement of (4.4) is obtained for free,
the identifications just mentioned, in turn, are {\it not\/} obtained for free,
as we have shown in this paper.

\smallskip\noindent
{\smc Remark 2.}
For a based simplicial set
$(X,*)$, the realization
$|CX|$ of the cone
$CX$ is naturally homeomorphic
to the realization 
$|X \wedge \Delta[1]|$
of
$X \wedge \Delta[1]$.
In fact,
a suitable subdivision
of $|CX|$ 
yields a realization of
$X \wedge \Delta[1]$.
It is tempting trying to
construct a homeomorphism between
$|DENK|$ and $|WK|$
in a combinatorial way
by inductively constructing the requisite maps
between the 
realizations of the
constituents of (4.1) and of the
corresponding terms in (2.1)
but we did not succeed in so doing.
The problem
is that
the realization
of the simplicial
monoid
$\Delta[1]$
does
not 
yield the geometric monoid structure on the interval $I$
coming into play in Section 1 above whence
the realization of an action $X \times \Delta[1] \to X$
of $\Delta[1]$
on a simplicial set $X$
is {\it not\/}
an $I$-action
on the realization of $X$
in the sense of
Section 1.
Rather, the realization
of the simplicial
monoid structure on
$\Delta[1]$
yields
the function from $I \times I$ to $I$ which sends
$(a,b)$ to $\min (a,b)$.
A suitable homeomorphism
identifies this monoid structure
with the one considered in Section 1 above.
For example, as pointed out by the referee,
one could take 
the function which
assigns
$(a \max (a,b),b\max(a,b)) \in I^2$
to
$(a,b) \in I^2$.
Further, 
the monoid structure
arising from the function $\min$
also gives rise to a monad
in the category of spaces
and with reference to it, the construction
(1.7)
can still
be carried out;
formally the same argument as that for
the proof of our main result
then identifies the limit (say) $LK$ of the resulting sequence of spaces with
the realization
$|DNK|$ of $DNK$
and, by virtue of (4.4),
$LK$ is naturally homeomorphic to $B|K|$.
However we do not see how this homeomorphism 
may be obtained
directly
since
we are unable to identify the
monad
in the category of spaces arising
from
the unit interval having the usual multiplication
as monoid structure with the other monad arising
from
the function $\min$ as monoid structure.
\bigskip
\centerline{References}
\medskip
\widestnumber\key{999}

\ref \no \cbergthe
\by C. Berger
\paper Une version effective du th\'eor\`eme
de Hurewicz
\paperinfo Th\`ese de doctorat, Universit\'e de Grenoble, 1991 
\endref

\ref \no \bottone
\by R. Bott
\paper On the Chern-Weil homomorphism and the continuous cohomology of
Lie groups
\jour Advances
\vol 11
\yr 1973
\pages  289--303
\endref

\ref \no \botshust
\by R. Bott, H. Shulman, and J. D. Stasheff
\paper On the de Rham theory of certain classifying spaces
\jour Advances
\vol 20
\yr 1976
\pages 43--56
\endref

\ref \no \curtione
\by E. B. Curtis
\paper Simplicial homotopy theory
\jour Advances in Math.
\vol 6
\yr 1971
\pages 107--209
\endref

\ref \no \dwykanon
\by W. G. Dwyer and D. M. Kan
\paper Homotopy theory and simplicial groupoids
\jour Indag. Math.
\vol 46
\yr 1984
\pages 379--385
\endref

\ref \no \eilmactw
\by S. Eilenberg and S. Mac Lane
\paper On the groups ${\roman H(\pi,n)}$. I.
\jour Ann. of Math.
\vol 58
\yr 1953
\pages  55--106
\endref
\ref \no \eilmacth
\by S. Eilenberg and S. Mac Lane
\paper On the groups ${\roman H(\pi,n)}$. II. Methods of computation
\jour Ann. of Math.
\vol 60
\yr 1954
\pages  49--139
\endref

\ref \no \gugenhtw
\by V.K.A.M. Gugenheim
\paper On the chain complex of a fibration
\jour Illinois J. of Mathematics
\vol 16
\yr 1972
\pages 398--414
\endref

\ref \no \gugenmay
\by V.K.A.M. Gugenheim and J.P. May
\paper On the theory and applications of differential
torsion products
\jour Memoirs of the Amer. Math. Soc.
\vol 142
\yr 1974
\endref

\ref \no \huebkade
\by J. Huebschmann and T. Kadeishvili
\paper Small models for chain algebras
\jour Math. Z.
\vol 207
\yr 1991
\pages 245--280
\endref

\ref \no \lattice
\by J. Huebschmann
\paper 
Extended moduli spaces, the Kan construction,
and lattice gauge theory
\jour Topology 
\vol 38 
\yr 1999 
\pages
555--596
\paperinfo {\tt dg-ga/9505005, dg-ga/9506006}
\endref

\ref \no \kanone
\by D. M. Kan
\paper On homotopy theory and c.s.s. groups 
\jour Ann. of Math. 
\vol 68 
\yr 1958 
\pages 38--53
\endref

\ref \no \maclafou
\by S. Mac Lane
\paper Milgram's classifying space as a tensor product of functors
\paperinfo in: 
The Steenrod algebra and its applications, F. P. Peterson, ed.
\jour Lecture Notes in Mathematics
\vol 168
\yr 1970
\pages 135--152
\publ Springer-Verlag
\publaddr Berlin $\cdot$ Heidelberg $\cdot$ New York
\endref

\ref \no \maclbotw
\by S. Mac Lane
\book Categories for the Working Mathematician
\bookinfo Graduate Texts in Mathematics 
\vol 5
\publ Springer
\publaddr Berlin $\cdot$ G\"ottingen $\cdot$ Heidelberg
\yr 1971
\endref

\ref \no \milgrone
\by J. Milgram
\paper The bar construction and abelian H-spaces
\jour Illinois J. of Math. 
\vol 11
\yr 1967
\pages 242-250
\endref

\ref \no \milnoru
\by J. Milnor
\paper Construction of universal bundles. I. II
\jour Ann. of Math.
\vol 63
\yr 1956
\pages 272--284, 430--436
\endref
\ref \no \milnothr
\by J. Milnor
\paper The realization of a semi-simplicial complex
\jour Ann. of Math.
\vol 65
\yr 1957
\pages 357--362
\endref

\ref \no \mooreone
\by J. Moore
\paper Comparison de la bar construction \`a la construction $W$
et aux complexes $K(\pi,n)$
\paperinfo Expos\'e 13
\jour S\'eminaire H. Cartan
\yr 1954/55
\pages 242-250
\endref

\ref \no \puppeone
\by D. Puppe
\paper Homotopie und Homologie
in abelschen Gruppen und Monoidkomplexen. I. II
\jour Math. Z.
\vol 68
\yr 1958
\pages 367--406, 407--421
\endref

\ref \no \quillone
\by D. Quillen
\paper Higher algebraic K-theory, I
\paperinfo in: Algebraic K-theory I, Higher K-theories,
ed. H. Bass,
Lecture Notes in Mathematics, No. 341
\publ Springer
\publaddr Berlin $\cdot$ Heidelberg $\cdot$ New York $\cdot$ Tokyo
\yr 1973
\pages 85--147
\endref

\ref \no \gsegatwo
\by G. B. Segal
\paper Classifying spaces and spectral sequences
\jour Publ. Math. I. H. E. S. 
\vol 34
\yr 1968
\pages 105--112
\endref

\ref \no \gsegathr
\by G. B. Segal
\paper Categories and cohomology theories
\jour Topology
\vol 13
\yr 1974
\pages 293--312
\endref

\ref \no\stashone
\by J. D. Stasheff
\paper Homotopy associativity of H-spaces.I
\jour Trans. Amer. Math. Soc.
\vol 108
\yr 1963
\pages 275--292
\moreref
\paper II
\jour Trans. Amer. Math. Soc.
\vol 108
\yr 1963
\pages 293--312
\endref

\ref \no \stasheig
\by J. D. Stasheff
\paper H-spaces and classifying spaces: Foundations and recent developments
\jour Proc. Symp. Pure Math.
\vol 22
\yr 1971
\pages 247--272
\publ American Math. Soc.
\publaddr Providence, R. I.
\endref

\ref \no \steenone
\by N. E. Steenrod
\paper Milgram's classifying space of a topological group
\jour Topology 
\vol 7
\yr 1968
\pages 349--368
\endref

\ref \no \weingone
\by S. Weingram
\paper The realization of a semisimplicial bundle map
is a $k$-bundle map
\jour Trans. Amer. Math. Soc.
\vol 127
\yr 1967
\pages 495--514
\endref

\ref \no \wongone
\by S.-C. Wong
\paper Comparison between the reduced bar construction and 
the reduced $W$-construction
\paperinfo Diplomarbeit, Math. Institut der Universit\"at Heidelberg,
1985
\endref

\enddocument